\documentclass[11pt,a4paper]{article}
\usepackage[english]{babel}
\usepackage{hyperref}
\usepackage[utf8]{inputenc}
\usepackage{color}

\usepackage{fontenc}

\usepackage{amssymb, amsmath, amsbsy}
\usepackage{array}
\usepackage{amsthm}
\usepackage{fancyhdr}
\usepackage{nccmath}
\usepackage{geometry}
\usepackage{graphicx} % figuras
\usepackage{subfigure} % subfiguras
\usepackage{float} % para que ponga la figura exactamente en el lugar en el que se da la sentencia, con [H]

\newcommand{\<}{\left<}
\renewcommand{\>}{\right>}

\newcommand{\com}[1]{\opt{draft}{\textcolor{red}{
$\LHD$ #1 $\RHD$\marginpar{\textcolor{red}{$\begin{lema}acksquare$}}}}}

\newcommand{\comb}[1]{\opt{draft}{\textcolor{blue}{
$\LHD$ #1 $\RHD$\marginpar{\textcolor{blue}{$\begin{lema}acksquare$}}}}}

\newenvironment{demo}{{\bf Proof }}
{\qed \\}

\newcommand{\re}{\mathbb R}

\renewcommand{\(}{\left(}
\newcommand{\lb}{\label}
\newcommand{\nn}{\nonumber}

\newcommand{\fracc}{\displaystyle\frac}
\newcommand{\ds}{\displaystyle}
\renewcommand{\)}{\right)}

\newcommand{\eps}{\ensuremath{\varepsilon}}

\def\p{\varphi}

\def\vec{\overrightarrow}

\def\k{\frak k}

\def\s{{\rm s_\lambda}}
\def\c{{\rm c_\lambda}}

\def\co{{\rm co_\lambda}}

\def\ta{{\rm ta_\lambda}}
\def\co{{\rm co_\lambda}}

\newcommand{\bde}{\begin{defi}}
\newcommand{\ede}{\end{defi}}

\newcommand{\be}{\begin{enumerate}}
\newcommand{\ee}{\end{enumerate}}

\newcommand{\ba}{\begin{array}}
\newcommand{\ea}{\end{array}}

\def\oMu{{\overline{M}_\lambda^{n+1}}}
\def\oMd{{\overline{M}_\lambda^{2}}}

\def\hM{{\widehat M }}

\def\hB{{\widehat {\rm B}}}

\def\hA{{\widehat {A}}}
\def\hC{{\widehat {C}}}

\newtheorem{defi}{\hspace{12pt} Definition}
\newtheorem{teor}{\hspace{12pt} Theorem}
\newtheorem{prop}[teor]{\hspace{12pt} Proposition}

\newtheorem{lema}[teor]{\hspace{12pt} Lemma}
\newtheorem*{lema*}{\hspace{12pt} Lemma}
\newtheorem*{teor*}{\hspace{12pt} Theorem}

\def\vec{\overrightarrow}

\def\k{\frak k}

\def\s{{\rm s_\lambda}}
\def\c{{\rm c_\lambda}}

\def\co{{\rm co_\lambda}}

\def\ta{{\rm ta_\lambda}}
\def\co{{\rm co_\lambda}}

\def\h{{\kappa}}

\newcommand{\ben}{\begin{enumerate}}
\newcommand{\een}{\end{enumerate}}
\newcommand{\bi}{\begin{itemize}}
\newcommand{\ei}{\end{itemize}}
\newcommand{\bec}{\begin{equation}}
\newcommand{\eec}{\end{equation}}
\newcommand{\beca}{\begin{equation*}}
\newcommand{\eeca}{\end{equation*}}
\newcommand{\bal}{\begin{align}}
\newcommand{\aal}{\end{align}}
\newcommand{\bala}{\begin{align*}}
\newcommand{\aala}{\end{align*}}

\begin{document}

\title{A discrete Blaschke Theorem for convex polygons in $2$-dimensional space forms \footnotetext{The research is partially supported by grant PID2019-105019GB-C21 funded by MCIN/AEI/ 10.13039/501100011033 and by \lq\lq ERDF A way of makimg Europe\rq\rq , and by the grant AICO 2021 21/378.01/1 funded by the Generalitat Valenciana.}
}
\author{Alexander Borisenko and Vicente Miquel}

\date{}

\maketitle

\begin{abstract}
Let $M$ be a $2$-space form. Let $P$ be a convex polygon in $M$. For these polygons, we define (and justify) a curvature $\h_i$ at each vertex $A_i$ of the polygon and and prove the following Blaschke's type theorem: \lq\lq If $P$ is a convex plygon in $M$ with curvature at its vertices $\h_i\ge \h_0 >0$, then the circumradius $R$ of $P$ satisfies $\ta(R) \le \pi/(2\h_0)$ and the equality holds if and only if the polygon is a $2$-covered segment\rq\rq
\end{abstract}

%\maketitle

\section{Introduction and main result}

We start recalling that a $n$-dimensional space form $\oMu$ of curvature $\lambda$ is a simply connected $n$-dimensional Riemannian manifold of constant sectional curvature $\lambda$. The only ones are: when $\lambda =0$, the Euclidean Space, when $\lambda>0$, the $n$-dimensional sphere of radius $1/\sqrt{\lambda}$ in the euclidean $\re^{n+1}$, and, when $\lambda<0$, the Hyperbolic Space of sectional curvature $\lambda$, that can be visualized as the upper connected component of the Minkowski sphere of radius $1/\sqrt{|\lambda|}$ in the Minkowski Space $\re^{n+1}$.

In the book of Blaschke \cite{Bl49} it is proved that
{\it  if $\Gamma$ is a closed convex regular curve in the Euclidean plane that bounds a compact convex region $\Omega$ and the curvature $\h$ of $\Gamma$ is bounded from below by some constant $\h_0>0$, then, for every point $p\in \Gamma$, the circle tangent to $\Gamma$ at $p$ and with radius $R=\fracc1{\h_0}$ bounds a disk that contains $\Omega$.}

This result was extended by H. Karcher (\cite{Ka68}) for the other space forms. Before stating it, we recall a notation that allows to describe the 
 geometry of  space forms can be described in a unified way:
\bec \lb{tri_hyp}
\s(t)= \begin{cases} \frac{\sin(\sqrt{\lambda}
t)}{\sqrt{\lambda}},  & \text{ if } \lambda>0 \\
t &  \text{ if }  \lambda =0 \\
\frac{\sinh(\sqrt{|\lambda|}
t)}{\sqrt{|\lambda|}} &  \text{ if } \lambda < 0  
\end{cases} , \quad
\ba{c} \c(t)=\begin{cases} \cos(\sqrt{\lambda}
t),  & \text{ if } \lambda>0 \\
1 &  \text{ if }  \lambda =0 \\
\cosh(\sqrt{|\lambda|}
t) &  \text{ if } \lambda < 0  
\end{cases},\\
\ta(t) = \ds \frac{\s(t)}{\c(t)},  \quad 
\text{and} \quad
\co(t) =
\ds\frac{\c(t)}{\s(t)}\ea.
\eec
The functions above satisfy the following computational rules:
\bec\label{trirel}
{\rm c}'_\lambda = - \lambda \, \s, \quad \s'(t) = \c(t),\quad {\rm c}^2_\lambda + \lambda\, {\rm s}^2_\lambda =1, \quad 
  \quad \frac{1}{\c^2(t)} = 1 + \lambda \ta^2(t),
  \eec
\bec \s(t+u) = \s(t) \c(u) + \c(t) \s(u) \text{ and }  \c(t+u) = \c(t) \c(u) - \lambda\ \s(t) \s(u).\nn
\eec
where \lq\lq\ $ ' $ \rq\rq denotes the derivative respect to $t$.

We shall recall also the following concept:

{\it Given any convex closed curve $\Gamma$ in $\oMd$,  the  \underline{circumradius} of $\Gamma$ is the minimum value of $R$ such that a disk of radius $R$ in $\oMd$ contains the domain bounded by $\Gamma$}. 

With this concept, the Blaschke-Karcher theorem can be stated in the following form:

\begin{teor}[\cite{Ka68}]\lb{BlTh}
if $\Gamma$ is a closed convex $C^2$ curve in $\oMd$  that bounds a compact convex region $\Omega$ and with curvature $\h$ (where it is $C^2$) satisfies $\h\ge\h_0>\max\{0,-\lambda\}$, then the circumradius $R$ of $\Gamma$ satisfies $\ta(R)\le\fracc1{\h_0}$. 
\end{teor}

Further developments of related Blaschke theorems has been done by obtaining conditions under which a convex set in $\re^n$ can be included  in other (\cite{Ra74,De79,BS89}) or its generalization to Riemannian manifolds where the included convex is a sphere (\cite{Ho99}).

In Theorem \ref{BlTh} (and in the other cited developments), the hypothesis of strong convexity ($\h\ge\h_0>0$) is necessary, the theorem is not true for $\h\ge 0$. Then  it cannot be applied to closed convex polygons. Here we shall show that it is possible to have a version of the theorem for polygons once we give an appropriate definition of curvature at the vertices of a polygon. We shall take the following one:

\begin{defi}\lb{defkA} Let $A$ be a vertex of a convex polygon $P$ in a space form $\oMd$. When $\lambda>0$, the sides $\ell_i$ of $P$ must satisfy $\ell_i < \pi/\sqrt{\lambda}$. Let $\hA$ be the interior angle of $P$ at the vertex $A$, and let $\ell_1,\ell_2$ be the lengths of the sides of $P$ that meet at vertex $A$. We define the \lq\lq curvature of $P$ at $A$\rq\rq  by the number
\bec \lb{kA}
\h_A = \frac{(\pi-\hA)}{\ta(\ell_1/2)+\ta(\ell_2/2)}. 
\eec
A compact convex polygon $P$ is called $\h_0$-convex if at every vertex $A$ of $P$, $\h_A \ge \h_0>0$.
\end{defi}
\bec \hskip-4.5cm\text{When $\lambda=0$, the definition \eqref{kA} becomes $\h_A = \fracc{2 (\pi-\hA)}{\ell_1+\ell_2}$.} \lb{kAE}
\eec

In the next section we shall give the reasons why we have chosen Definition \eqref{kA}. 

The version of Theorem \ref{BlTh} that shall prove for polygons is:

\begin{teor}\lb{BTE}
Let $P$ be a compact $\h_0$-convex polygon in $\oMd$, with sides lower than $\pi/\sqrt{\lambda}$ if $\lambda>0$, and with curvature at each vertex $A_i$ satisfying $\h_{A_i} \ge\h_0 > \max\{0,-\lambda\}$. Then the circumradius $R$ of $P$  satisfies 
\bec\lb{inTh}\ta(R) \le \pi/(2\h_0),
\eec
 and the equality  holds if and only if the polygon degenerates to a  $2$-covered segment.
\end{teor}

For the euclidean plane, Definition \eqref{kAE} was used in \cite{BCM} in the study of approximations of surfaces by planar triangulations and in \cite{BO10} it was studied how good is this definition to approximate the curvature of a planar curve by a polygonal line. Other applications of this definition in the euclidean case has been done in \cite{CRR}. Related but different definitions of curvature of a polygon in the euclidean plane have been used for other applications in \cite{Mu04} and \cite{SM}.

Some people could prefer to take \eqref{kAE} as definition for the curvature of a convex polygon for every $\oMd$, without having into account the value of $\lambda$. In the last section of the paper we shall give the corresponding result (Th. \ref{BTEb}) for this definition.

\section{About Definition \ref{defkA}}

If we consider a convex polygon as a limit of smooth curves approaching it and the curvature at a vertex as the limit of the curvature of the points at the curves whose limit is the vertex, then the curvature becomes infinite. Obviously, this is not a good definition for many geometric properties. We would like a definition satisfying the following properties:

\begin{itemize}
\item[P1] the curvature of a vertex is bigger as the interior angle is lower
\item[P2] the curvature of a vertex is bigger as the lengths of the adjacent sides is shorter
\item[P3] if we have a regular polygon inscribed in a circle and we take the number of sides of the regular polygon increasing up to infinite, the curvature of the vertices approach the curvature of the circle.
\end{itemize} 

Properties P1 and P3 correspond to a natural geometric intuition. Property P2 is related to the fact that we want to generalize Theorem \ref{BlTh} which fails when $\h_0=0$ (in the Euclidean case) because with $\h_0=0$ you may have straight lines with arbitrary length which are the obstacle for upper bounds for the circumradius.

It is obvious that our definition \ref{defkA} satisfies P1 and P2. In the next Proposition, we shall check that it also satisfies P3.

\begin{prop}\lb{Pro}
Let $C$ be a circle of radius $R$ ($<\pi/(2\sqrt{\lambda})$ if $\lambda>0$) in $\oMd$ and let $P_n$ be a regular polygon of $n$ sides with its vertices in $C$. If $\h_n$ denotes the specific curvature at any vertex $A_n$ of the polygon $P_n$, then $\lim_{n\to \infty} \h_n = \co(R)$, which is the curvature of $C$. 
\end{prop}
\begin{demo}
Let us recall some trigonometric formulae of the space forms: Let $\Delta$ be a geodesic triangle with sides $a,b,c$ and opposite vertices $A,B,C$. Let $\hA,\hB,\hC$ be the angles at these vertices. Then the following formulae hold:
\begin{align}
\cos \hA&= \frac{\c(a) - \c(b)\ \c(c)}{\lambda\ \s(b)\ \s(c)} \lb{ca}\\
(\text{when } & \lambda=0 \text{ the quotient in the second term of \eqref{ca} must be understood}\nn \\ \text{taking} &\text{ limits for } \lambda\to 0, \text{giving the standard cosinus law in Euclidean plane}) \nn \\
\frac{\sin\hA}{\s(a)} &= \frac{\sin\hB}{\s(b)} = \frac{\sin\hC}{\s(c)} \lb{sinl} %\\
%\cos \hA &= - \cos\hB\ \cos \hC + \sin\hB\ \sin\hC\ \c(a) \lb{cosA}
\end{align}
Let $A_n$, $B_n$ be two consecutive vertices of the polygon, bounding a side $A_nB_n$ of length $\ell_n$. Let $M_n$ be the middle point of $A_nB_n$ between $A_n$ and $B_n$. Let $O$ be the center of the circle. Let us consider the geodesic triangle $OM_nB_n$, and denote by $d_n$ the length of the geodesic $OM_n$. One has $\hM_n=\pi/2$. We can apply \eqref{ca} and \eqref{sinl} to this triangle to obtain:
\begin{align*}
\c(R)&= \c(d_n) \c(\ell_n/2), \qquad 
\s(d_n) = \s(R) \sin(\hB_n/2)
\end{align*} 
From these two equalities and the formulae \eqref{trirel} we obtain
\begin{align}
&\frac{\c^2(R)}{\c^2(\ell_n/2)} = \c^2(d_n) = 1-\lambda \s^2(d_n) = 1-\lambda\ \s^2(R) \sin^2(\hB_n/2)\text{, then } \nn\\
&\sin^2(\hB_n/2) =  \frac1{\lambda\ \s^2(R)}-\frac{\co^2(R)}{\lambda\ \c^2(\ell_n/2)} = 1 + \frac{\co^2(R)}{\lambda} \(1-\frac1{\c^2(\ell_n/2)}\) \nn\\
&\cos^2((\hB_n/2) = 1- \sin^2(\hB_n/2) = \co^2(R) \ta^2(\ell_n/2), \text{ then} \nn \\
\lb{ellbeta}
&\ta(\ell_n/2) = \ta(R) \cos(\hB_n/2) = \ta(R) \sin(\pi/2-\hB_n/2)
\end{align} 
We apply now the definition \eqref{kA} to the curvature $\k_n$ of $\hB_n$
\bec \lb{kn}
\h_n = \frac{(\pi-\hB_n)}{2 \ta(\ell_n/2)} = \frac{2(\pi/2-\hB_n/2)}{2 \ta(R) \sin(\pi/2-\hB_n/2)}. 
\eec
But $\lim_{n\to\infty} \hB_n = \pi$, then the limit of the quotient in \eqref{kn} for $n\to\infty$ is $\co(R)$, as claimed in Prop. \ref{Pro}.  \end{demo}

Definition \ref{defkA} is not the unique satisfying properties P1 to P3. If we take \eqref{kAE} as a definition of the curvature at $A$ for any value of $\lambda$, it is obvious that it satisfies P1 and P2, moreover P3 follows with the same proof of Proposition \ref{Pro} , having into account that $\lim_{n\to\infty}\frac{\ta(\ell_n/2)}{(\ell_n/2)} =1$. We have preferred \eqref{kA} because it gives a clean bound \eqref{inTh} in Theorem \ref{BTE}, but some people could prefer the other definition.

\section{Proof of the Theorem}

Let $A_1, A_2, ..., A_n$ be the consecutive vertices of the polygon $P$. Let $\hat A_i$ be the angles at these vertices, and $\ell_1, ..., \ell_i, ... , \ell_n$ be the lengths of the sides $A_nA_1, ... , A_{i-1}A_i, ... , A_{n-1} A_n$ respectively (see Figure 1).

For every segment $A_{i-1}A_i$ we construct a segment of circle $C_i$ of radius $\rho_i$ with center $O_i$ in a line orthogonal to $A_{i-1}A_i$ in its middle point and in the ray in the inward direction, and with boundary points $A_i$ and $A_{i+1}$. 

\begin{figure}[!h]
\begin{center}
 \includegraphics[scale=0.6]{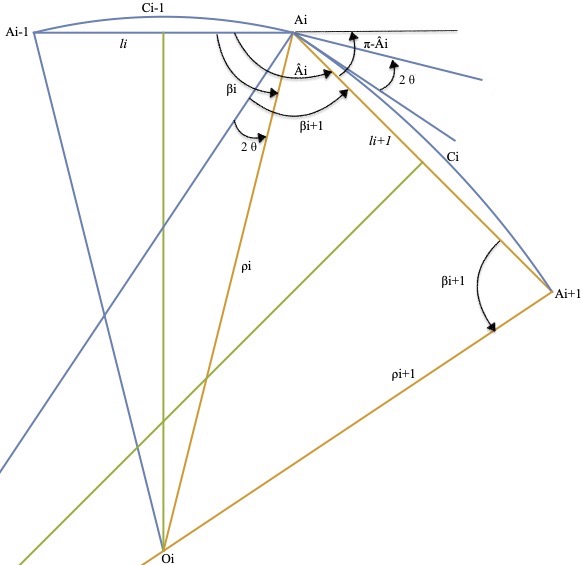}
 \caption{}\lb{fig1}
 \end{center}
\end{figure}

Each angle $\beta_i$ at $A_i$ of the isosceles triangle $A_{i-1}A_iO_i$ satisfies $0<\beta_i<\pi/2$ and the analogous of \eqref{ellbeta} for this triangle is 
\bec\lb{ellbeta2}
\ta(\ell_i/2) = \ta(\rho_i)\cos \beta_i = \ta(\rho_i)\sin(\delta_i)
\eec
 if we take $\delta_i = \pi/2-\beta_i$.

We now take the curve obtained as the union of the segments of circle $C_i$. This curve will be convex if and only if, for every $i=1, .., n$, the tangent vectors at $A_i$ of the circles $C_{i-1}$ and $C_i$ with the curve $C_{i-1}$ oriented from $A_{i-1}$ to $A_i$ and the curve $C_i$ oriented from $A_{i}$ to $A_{i+1}$ form an angle $2\ \theta_i $ in the interval $[0,\pi]$. This angle is the same that the one formed by the normals at $A_i$ to $C_{i-1}$ and $C_i$ pointing inward. These normals are $A_iO_i$ and $A_iO_{i+1}$, and this angle is non negative if and only if $\beta_i+\beta_{i+1} \ge  \hat A_i$, that is $\pi-\delta_i-\delta_{i+1}\ge  \pi-(\pi-\hat A_i)$, 
\bec\lb{diA}
\delta_i+\delta_{i+1} \le 
\pi-\hat A_i.
\eec

For every $i$, let us choose $\rho_i=R$ such that $\ta(R) = \frac{c}{\h_0}$.  From the hypothesis of $\h_0$-convexity of $P$, and using formula \eqref{ellbeta2}, we have 
\begin{align}
\h_0 \le & \frac{\pi-\hat A_i}{\ta(\ell_{i}/2)+\ta(\ell_{i+1}/2)} = \frac{ (\pi -\hat A_i)}{\ta(\rho_i) \sin \delta_i+ \ta(\rho_{i+1}) \sin \delta_{i+1}} = \frac{\h_0}{c}\  \frac{ (\pi-\hat A_i)}{\sin \delta_i+ \sin \delta_{i+1}},  \nn \\
\text{ then } & \pi - \hat A_i \ge c\ (\sin \delta_i+ \sin \delta_{i+1}) \ge c \frac2\pi \ (\delta_i + \delta_{i+1}) \lb{bpAi}
\end{align}
an inequality which coincides with \eqref{diA} if $c=\pi/2$, which gives $\ta(\rho_i) = \ta(R)=\fracc{\pi}{2\  \h_0} $. 
As a consequence the closed curve $C$ formed by the union of the $C_i$ is convex and with curvature equal $\fracc{2 \h_0}{\pi}$ at every regular point.  Now, for every $i$, let us take the exterior parallel curve $C_{i\eps}$ at distance $\eps$ from $C_i$. Each $C_{i\eps}$ is an arc of circle of radius $R+\eps$, then it has curvature  
$$\h_i = \co(R+\eps)=\fracc{\co(R)+\lambda \ta(\eps)}{1-\co(R) \ta(\eps)} = \fracc{2 \h_0/\pi+\lambda \ta(\eps)}{1-\ta(\eps) 2 \h_0/\pi }.$$

\begin{figure}[!h]
\begin{center}
 \includegraphics[scale=0.75]{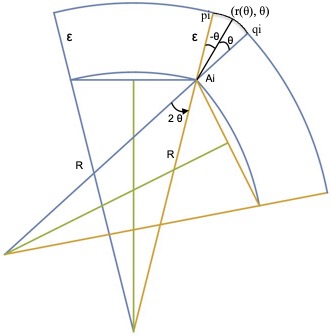}
 \caption{}\lb{fig1}
 \end{center}
 \end{figure}
 
Given to consecutive arcs $C_{\eps i-1}$ and $C_{\eps i}$, we join the right extreme $p_i$ of $C_{\eps i-1}$ with the left extreme $q_i$ of $C_{\eps i}$ by a curve that, in polar geodesic coordinates $(r,\p)$ around $A_i$ has the expression $c_i(\p)=(r_i(\p),\p)$, with $-\theta_i\le \p\le\theta_i$ (see Figure 2), satisfying

\begin{align}
p_i &= c_i(-\theta_i) = (r_i(-\theta_i),-\theta_i), \quad q_i = c_i(\theta_i) = (r_i(\theta_i),\theta_i) \nn\\ &\text{ that is } \quad r_i(-\theta_i)=\eps = r_i(\theta_i) \lb{c1}\\
c_i'(-\theta_i) &\text{ is tangent to } C_{(i-1)\eps} \text{ at } p_i \text{ and } c_i'(\theta_i) \text{ is tangent to } C_{i\eps} \text{ at } q_i \nn \\
&\text{that is } r_i'(-\theta_i)= 0 = r_i'(\theta_i) \lb{c2}\\
\text{the curvature $k_i$ of }& c_i \text{ at } c_i(-\theta_i) \text{ is equal to the curvature of } C_{(i-1)\eps} \text{ at } p_i \text{ and }\nn \\
\text{ the curvature $k_i$ of }& c_i \text{ at } c_i(\theta_i) \text{ is equal to the curvature of } C_{i\eps} \text{ at } p_i \nn \\
&\text{that is } k_i(-\theta_i)= \h_i =\co(R+r) = k_i(\theta_i). \lb{c3}
\end{align}

With these conditions, the curve $C_\eps$ obtained joining all the curves $C_{i\eps}$ and $c_i$ is $C^2$. Moreover, the function $r_i(\p)$ can be taken (see the appendix) such that, 
\bec
\lim_{\eps\to 0^+} k_i(\p)= \infty \quad \text {for } -\theta_i < \p < \theta_i,\lb{c4}
\eec
then for small enough $\eps$, $k_i>\co(R+\eps)$, then the curve $C_\eps$ has curvature bigger than $2\ \h_0/\pi + \psi(\eps)$ with $\psi(\eps) = \fracc{2 \h_0 + \lambda \pi^2}{\pi^2-\pi\ta(\eps) 2 \h_0} \ta(\eps),$ then $\lim_{\eps\to0}\psi(\eps)= 0$. 

 To each one of these curves $C_\eps$ we can apply Theorem \ref{BlTh} and conclude that the circumradius $R_\eps$ of each one of them satisfies $\ta(R_\eps) \le \fracc1{2\ \h_0/\pi + \psi(\eps)}$. Then, taking $\eps\to 0$, we obtain the inequality \eqref{inTh} of Theorem \ref{BTE}.

The equality holds if and only if equalities hold in all the inequalities of the above argument. In particular, equality implies $\sin \delta_i = \frac2\pi \delta_i$, which happens if and only if $\delta_i=\pi/2$. The other equalities that we must have are $\pi - \hat A_i = c\ (\sin \delta_i+ \sin \delta_{i+1}) = \pi$, that is $A_i=0$, and $\ta(\ell_i/2) = \ta(\rho_i) \sin \delta_i = \pi/\h_0$ and this is satisfied only in a $2$-covered segment of length $\pi/\h_0$. It is a degenerate polygon of curvature $\pi/(2 \ta(\ell_i/2)) = \h_0$.   \qed

\section{If we adopt definition \eqref{kAE}}

If the curvature at a vertex of a convex polygon is defined by \eqref{kAE}, then Theorem \ref{BTE} must be changed by:

\begin{teor}\lb{BTEb}
Let $P$ be a compact $\h_0$-convex polygon in $\oMd$ such that, if $\lambda>0$, the sides satisfy $\ell_i \le 2\ \frak e <\pi/\sqrt{\lambda}$ and, if $\lambda<0$, one has $\h_0 > \sqrt{-\lambda}$ and $\ell_i \ge 2\ \frak e >0$. Then the circumradius $R$ of $P$  satisfies 
\begin{align}
\lb{inThb1} R \le \pi/(2\h_0), \text{ if } \lambda=0\\
\lb{inThb2} \ta(R) \le \fracc{\ta(\frak e)}{\frak e} \fracc\pi{2\ \h_0} , \text{ if } \lambda \ne 0
\end{align}
 and the equality  holds if and only if the polygon degenerates onto a  $2$-covered segment.
 
 Let us observe that, when $\lambda<0$, $\ds \fracc{\ta(\frak e)}{\frak e} \fracc\pi{2\ \h_0} < \pi/(2\h_0)$ for $\frak e >0$.
\end{teor}

The proof is exactly the same for $\lambda=0$. For $\lambda >0$, one has that 
\bec
 \ell_i/2 = \frac{\ell_i/2}{\ta(\ell_i/2)}\ta(\ell_i/2) \ge \frac{\frak e}{\ta(\frak e)}\ta(\ell_i/2), \nn
 \eec
 then
 \begin{align}
&\h_0 \le  \frac{\pi-\hat A_i}{(\ell_{i}/2)+(\ell_{i+1}/2)} \le \frac{ (\pi -\hat A_i)}{\frac{\frak e}{\ta(\frak e)}\(\ta(\rho_i) \sin \delta_i+ \ta(\rho_{i+1}) \sin \delta_{i+1}\)}\nn \\
& \qquad = \frac{\h_0 \ta(\frak e)}{c\ \frak e}\  \frac{ (\pi-\hat A_i)}{\sin \delta_i+ \sin \delta_{i+1}},  \nn \\
&\text{ from which  }  \pi - \hat A_i \ge c\ \frac{\frak e}{\ta(\frak e)} (\sin \delta_i+ \sin \delta_{i+1}) \ge c\ \frac{\frak e}{\ta(\frak e)} \frac2\pi \ (\delta_i + \delta_{i+1}), \lb{bpAi}
\end{align}
and the union of the arcs $C_i$ is convex if we take $c= \fracc{\ta(\frak e)}{\frak e} \fracc\pi2$, and the rest of the proof follows as for Th. \ref{BTE}.

The proof for the case $\lambda<0$ follows the same steps, with the unique change that now the function $\fracc{\ell/2}{\ta(\ell/2)}$ is decreasing now and we have to bound taking the minimum value of $\ell/2$.

\section*{Appendix: About formulae \eqref{c1}, \eqref{c2}, \eqref{c3}, \eqref{c4}}

For the convenience of the reader, in this appendix we give details of the following:

1) Computation of the curvature $k_i$ of a curve $c_i(\p)=(r_i(\p),\p)$ expressed in polar geodesic coordinates $(r,\p)$ centered in $A_i$.

In these coordinates, the expression of the metric is $dr^2 + \s(r) d\p^2$, the tangent vector to the curve $c_i(\p)$ is $c_i'(\p) = r_i'(\p) \partial_r + \partial_\p$, then the unit tangent vector $\vec{t}$ and the unit normal vector $\vec{n}$ pointing inside are:
\begin{align*}
\vec{t} &= \frac{1}{\sqrt{r_i'(\p)^2 +\s(r_i(\p))^2}}  \(r_i'(\p) \partial_r + \partial_\p \), \\
 \vec{n} &= \frac{1}{\s(r_i(\p)) \sqrt{r_i'(\p)^2 +\s(r_i(\p))^2}} \(-\s^2(r_i(\p))\partial_r +r_i'(\p) \partial_\p \).
\end{align*}
Having into account that for the covariant derivative in $\oMd$ in geodesic polar coordinates we have 
\begin{align*}
\nabla_{\partial_r} \partial_\p &= \nabla_{\partial_\p} \partial _r = \co(r) \partial_\p, \quad \nabla_{\partial_r} \partial_r =0,\\
\nabla_{\partial_\p} \partial_\p & = \<\nabla_{\partial_\p} \partial_\p,{\partial_r}\> \partial_r + \< \nabla_{\partial_\p} \partial_\p ,\frac{1}{\s(r)} \partial_\p \> \frac{1}{\s(r)} \partial_\p = - \s \c \partial_r + \c r_i' \frac{1}{\s(r)} \partial_\p
\end{align*}
a straightforward computation gives
\begin{align}\lb{ki}
k_i(\p) = \<\nabla_{\vec{t}} \vec{t}, \vec{n}\> = \frac{\s}{(\s^2 + r_i'^2)^{3/2}} \( r_i'' + \s \c + 3 r_i'^2 \co\)
\end{align}

%\begin{center}
 %\includegraphics[scale=0.75]{parallelsa.jpg}
% \end{center}
 
2) Taking $r_i (\p)= r_{io} + a \p^2 + b\p^4$, determination of the coefficients  $r_{io} ,\  a ,\  b$ such that \eqref{c1}, \eqref{c2} and \eqref{c3} be satisfied.

%\begin{figure}[!h]
%\begin{center}
% \includegraphics[scale=0.75]{parallelsa.jpg}
% \caption{}\lb{fig1}
% \end{center}
%\end{figure}

With this definition it is obvious that $r_i$ satisfies $r_i(-\p)=r_i(\p)$, $r_i'(-\p)=- r_i'(\p)$ and $r_i''(-\p)=r_i''(\p)$. Conditions \eqref{c1} and  \eqref{c2} are satisfied if and only if
\begin{align}\lb{a}
 r_{io} = \eps - a \theta_i^2 - b \theta_i^4 \text{ and } a= -2 b \theta_i^2,\text{ then } r_{io}= \eps+ b\theta_i^4
 \end{align}
 From these it follows that $r_i''(\theta_i) = -8 b \theta_i^2$. Substitution of this in \eqref{ki} and application of equality \eqref{c3} gives:
 \begin{align}\lb{b}
 b &= \frac{-\s(\eps)}{8 \theta_i^2} \(\c(\eps) - \frac{\c(R) \c(\eps) \s(\eps) -\lambda \s(\eps)^2 \c(R)}{\c(\eps) \s(R) + \c(R) \s(\eps)} \)
 \end{align} 
3) Using the above values of $a$ and $b$ and the formula for $k_i$, checking that \eqref{c4} is satisfied. In fact, it follows from \eqref{b} that
\begin{align*}
\lim_{\eps\to 0^+} \frac{b}{\eps} &= \frac{-1}{8 \theta_i^2} , \qquad \lim_{\eps\to 0^+} \frac{a}{\eps} = \frac14 %\( 1 - \frac{\c(R)  }{ \s(R) }\),
, \qquad 
\lim_{\eps\to 0^+} \frac{r_{io}}{\eps} = 1 -\frac18 \theta_i^2\\
\lim_{\eps\to 0^+} \frac{r_i(\p)}{\eps} &= \(1 - \frac1{8}  \theta_i^2\) + \frac14 \p^2 - \frac1{8\theta_i^2} \p^4 = 1- \frac{1}{8\theta_i^2} (\theta_i^2-\p^2)^2\\
  \lim_{\eps\to 0^+} \frac{r_i'(\p)}{\eps} &= \(\frac{\p}2 + \frac{- \p^3}{2 \theta_i^2}\), \qquad
   \lim_{\eps\to 0^+} \frac{r_i''(\p)}{\eps} = \( \frac{1}2 + \frac{- 3 \p^2}{2 \theta_i^2} \)   \end{align*}
From \eqref{ki} and all the above limits,

\begin{align}
\lim_{\eps\to 0^+} k_i(\p)&= \lim_{\eps\to 0} \frac{\s}{(\s^2 + r_i'^2)^{3/2}} \( r_i'' + \s \c + 3 r_i'^2 \co\) \approx \nn \\
%%%%%%%%%
\approx & \lim_{\eps\to 0} \frac{\(1 - \frac1{8 \theta_i^2} (\theta_i^2-\p^2)^2 \)\eps }{\(\(1 - \frac1{8 \theta_i^2} (\theta_i^2-\p^2)^2 \)^2\eps^2  + \(\frac{\p}2 + \frac{- \p^3}{2 \theta_i^2}\)^2  \eps^2\)^{3/2}} \times\nn \\
%%%%%
&\qquad \times \( \( \frac{1}2 + \frac{- 3 \p^2}{2 \theta_i^2} \)  \eps + \(1 - \frac1{8 \theta_i^2} (\theta_i^2-\p^2)^2\) \eps %\right. \nn \\
%%%%
%&\left. \qquad \qquad \qquad  
+ 3 \frac{\(\frac{\p}2 + \frac{- \p^3}{2 \theta_i^2}\)^2  \eps^2}{\(1 - \frac1{8 \theta_i^2} (\theta_i^2-\p^2)^2\) \eps}\) \nn \\
%%%%%%%%%%
\approx & \lim_{\eps\to 0} \frac{\(1 - \frac1{8 \theta_i^2} (\theta_i^2-\p^2)^2 \) }{\(\(1 - \frac1{8 \theta_i^2} (\theta_i^2-\p^2)^2 \)^2  + \(\frac{\p}2 + \frac{- \p^3}{2 \theta_i^2}\)^2  \)^{3/2} \eps} \times\nn \\
%%%%%
& \qquad \times \( \( \frac{1}2 + \frac{- 3 \p^2}{2 \theta_i^2} \)   + \(1 - \frac1{8 \theta_i^2} (\theta_i^2-\p^2)^2\)  %\right. \nn \\
%%%%
%&\left. \qquad \qquad \qquad
  + 3 \frac{\(\frac{\p}2 + \frac{- \p^3}{2 \theta_i^2}\)^2 }{\(1 - \frac1{8 \theta_i^2} (\theta_i^2-\p^2)^2\)}\) \nn 
\end{align}
And, having into account that $\theta_i \in[0,\pi/2]$, \eqref{c4} follows easily for $\p \in]-\theta_i,\theta_i[$. 

\bibliographystyle{alpha}

\begin{thebibliography}{99}
 

\bibitem {Bl49}  Wilhelm  Blaschke; {\it Kreis und Kugel}. Chelsea Publishing Co., New York, 1949. x+169 pp. (Photo-offset reprint of the edition of 1916 [Veit, Leipzig]) 

\bibitem{BCM} Vincent Borrelli, Frederic Cazals and Jean Marie Morvan;  On the angular defect of triangulations and the pointwise approximation of curvatures. {\it Computer Aided Geometric
Design} 20, (2003) 319–341.

\bibitem{BO10} Vincent Borrelli and Fabrice Orgeret;  Error term in pontwise approximation of the curvature of a curve, {\it Computer Aided Geometric Design} 27 (2010) 538-550.

\bibitem{BS89} Jeff Brooks and John B. Strantzen; Blaschke’s rolling theorem in $\re^n$.
{\it Mem. Amer. Math. Soc.} 80 (1989), no. 405, vi+101 pp.

\bibitem{CRR} Juliá Cufí, Agustí Reventós and Carlos J. Rodríguez; Curvature for Polygons, {\it Amer. Math. Monthly},122,4 (2015) 332-337

\bibitem{De79} José A. Delgado; Blaschke's theorem for convex hypersurfaces. {\it J. Differential Geometry} 14 (1979), no. 4, 489–496 (1981). 

%\bibitem {BD13} Borisenko, A. A.; Drach, K. D. On the sphericity of hypersurfaces with normal curvature bounded from below. (Russian) {\it Mat. Sb.} 204 (2013), no. 11, 21–40; translation in {\it Sb. Math.} 204 (2013), no. 11-12, 1565–1583rst 2008. {\color{red} por ahora no hay ninguna referencia a este paper dentro del artículo}

\bibitem{Ho99} Ralph Howard; Blaschke rolling theorems for manifolds with boundary, {\it Manuscripta Math.} 99 (1999) 4, 471-483. 

\bibitem{Ka68}  Hermann Karcher; Umkreise und Inkreise konvexer Kurven in der sphärischen und der hyperbolischen Geometrie. {\it Math. Ann.} 177 (1968), 122–132.

\bibitem{Mu04} O. R. Musin; Curvature extrema and four-vertex theorems for polygons and polyhedra, {\it Journal of Mathematical Sciences}, 119,2 (2004) 268--277.

\bibitem{Ra74} Jeffrey Rauch; An inclusion theorem for ovaloids with comparable second fundamental forms. {\it J. Differential Geometry} 9 (1974), 501–505. 

\bibitem{SM} Koya Sakakibara  and Yuto Miyatake, Yuto; A fully discrete curve-shortening polygonal evolution law for moving boundary problems, {\it Journal of Computational Physics } 424 (2021) 109857 (22 pages).


\end{thebibliography}

{ Alexander Borisenko, \\
{B. Verkin Institute for Low Temperature, Physics and Engineering of the National Academy of Sciences of Ukraine\\
Kharkiv, Ukraine
}\\ and \\
 Department of Mathematics \\
University of Valencia \\
46100-Burjassot (Valencia), Spain}

{aborisenk@gmail.com}\\

{ Vicente Miquel \\
Department of Mathematics \\
University of Valencia \\
46100-Burjassot (Valencia), Spain
}

{vicente.f.miquel@uv.es}

\end{document}